\def\eps{\epsilon}
\def\tilde{\widetilde}

\baselineskip=14pt
\parskip=10pt
\def\Tilde{\char126\relax}
\def\halmos{\hbox{\vrule height0.15cm width0.01cm\vbox{\hrule height
 0.01cm width0.2cm \vskip0.15cm \hrule height 0.01cm width0.2cm}\vrule
 height0.15cm width 0.01cm}}
\font\eightrm=cmr8  
\font\eighttt=cmtt8
\magnification=\magstephalf

\parindent=0pt
\overfullrule=0in
\bf
\centerline
{PROOF OF THE REFINED ALTERNATING SIGN MATRIX CONJECTURE}
\rm
\bigskip
\centerline{ {\it Doron ZEILBERGER}\footnote{$^1$}
{\eightrm  \raggedright
Department of Mathematics, Temple University,
Philadelphia, PA 19122, USA. \break
E-mail: {\eighttt zeilberg@math.temple.edu} ;
WWW: {\eighttt http://www.math.temple.edu/\Tilde zeilberg}.
Supported in part by the NSF. Oct. 29,1995.
This version: March 20, 1996. To appear in New York J. of
Math.
} 
}
 
{\bf Abstract.} 
Mills, Robbins, and Rumsey conjectured, and Zeilberger proved, that
the number of alternating sign matrices of order $n$ equals
$A(n):={{1!4!7! \cdots (3n-2)!} \over {n!(n+1)! \cdots (2n-1)!}}$.
Mills, Robbins, and Rumsey also made the stronger conjecture that
the number of such matrices
whose (unique) `$1$' of the first row is at the $r^{th}$ column,
equals
$A(n) {{n+r-2} \choose {n-1}}{{2n-1-r} \choose {n-1}}/
{{3n-2} \choose {n-1}}$.
Standing on the shoulders of A.G. Izergin, V. E. Korepin, and G. Kuperberg,
and using in addition orthogonal polynomials  and $q$-calculus,
this stronger conjecture is proved.
 
{\bf INTRODUCTION}
 
An {\it alternating sign matrix}, or ASM, is a matrix of $0$'s, $1$'s,
and $-1$'s such that the non-zero elements in each row and each column
alternate between $1$ and $-1$ and begin and end with $1$, for example:
$$
\pmatrix{
0&1&0&0\cr
1&-1&1&0\cr
0&0&0&1\cr
0&1&0&0\cr} \quad .
$$
 
Mills, Robbins, and Rumsey
[MRR1][MRR2]([S], conj. 1) conjectured, and I proved[Z],
that there are
$$
A(n):={{1!4!7! \cdots (3n-2)!} \over {n!(n+1)! \cdots (2n-1)!}}
$$
alternating sign matrices of order $n$. Another, shorter, proof was
later given by Greg Kuperberg[K]. Kuperberg deduced the
straight enumeration of ASMs from their weighted enumeration
by Izergin and Korepin[KBI]. In this paper, I extend Kuperberg's
method of proof to prove the more general, refined enumeration, also
conjectured in [MRR1][MRR2], and listed by Richard Stanley[S]
as the third of of his ``Baker's Dozen'', that:
 
{\bf Main Theorem.}$-$ There are
$$
A(n,r):=A(n)  {{{{n+r-2} \choose {n-1}}{{2n-1-r} \choose {n-1}}} \over
{ {{3n-2} \choose {n-1}}}}
$$
$n \times n$ alternating sign matrices for which 
the (unique) `$1$' of the first row is at the $r^{th}$ column.  
 
As in Kuperberg's proof, we are reduced to evaluating a certain determinant.
Unlike the original determinant, it is not evaluable in closed form.
We evaluate it by using the q-analog of the Legendre polynomials
over an interval, 
introduced, and  greatly generalized, by Askey and Andrews[AA]
and Askey and Wilson[AW]. All that is needed from the general theory of
orthogonal polynomials and from $q$-calculus is reviewed, so 
a sufficient condition for following the present paper is 
having read Kuperberg's
paper[K]  that includes a very clear exposition of the Izergin-Korepin
formula, and of its proof. Of course, this is also a necessary condition.
In order to encourage readers to look up and read Kuperberg's beautiful
paper[K], and to save myself some typing,
I will use the notation, and results, of [K], without reviewing 
them.\footnote{$^2$}
{\eightrm Kuperberg's paper may be obtained from his Home Page
{\eighttt http://www.math.yale.edu/users/greg}.
Once you downloaded the file {\eighttt asm.ps.gz} you should first
get rid of the {\eighttt .gz} by typing {\eighttt gunzip asm.ps}, and then
to get a copy, you type {\eighttt lpr asm.ps}.}
 
{\bf Boiling It Down To a Determinant Identity}
 
Let $B(n,r)$ be the number of ASMs of order $n$ whose sole
`$1$' of the first (equivalently last) row, is at the $r^{th}$ column.
In order to stand on Kuperberg's shoulders more comfortably,
we will consider the last row rather than the first row.
 
Like in [K], define $[x]:=(q^{x/2}-q^{-x/2})/(q^{1/2}-q^{-1/2})$,
take $q:=e^{2 \sqrt{-1} \pi /3}$, and consider 
$$
Z(n\,;\,2, \dots ,\,\,2\,,\,2+a\,;\, 0, \dots , 0) \quad ,
$$
where between the two semi-colons inside $Z$ there are
$n-1$ $2$'s followed by a single $2+a$. Here $a$ is an indeterminate.
Let's look at an ASM of order $n$, whose sole `$1$' of the
last row is at the $r^{th}$ column. It is readily seen that
the $r-1$ zeros to the left of that `$1$' each contribute a
weight of $[2+a]$, while the remaining $n-r$ zeros, to the
right of the aforementioned `$1$', each contribute a weight of $[1+a]$.
The `$1$' itself contributes $q^{-1-a/2}$, which is $q^{-a/2}$ times
what it did before. Hence
$$
Z(n\,;\,2, \dots ,\,\,2\,,\,2+a\,;\, 0, \dots , 0)\,=\,
(-1)^n q^{-n} q^{-a/2} \sum_{r=1}^{n} B(n,r) [2+a]^{r-1}[1+a]^{n-r}
\quad .
$$
We also have, thanks to [K], (or [Z]:just plug in $a=0$ above):
$$
Z(n\,;\,2, \dots ,2\,,\,2\,;\, 0, \dots , 0)\,=\,
(-1)^n q^{-n} A(n) \quad .
$$
Hence:
$$
{{Z(n\,;\,2, \dots ,2\,,\,2+a\,;\, 0, \dots , 0)} \over
{Z(n;2, \dots ,2,\,2\,; 0, \dots , 0)}}=
{{q^{-a/2} } \over {A(n)}} 
\sum_{r=1}^{n} B(n,r) [2+a]^{r-1}[1+a]^{n-r} \quad.
$$
Since 
$$
\{ [2+a]^{r-1}[1+a]^{n-r} \,;\, 1 \leq r \leq n \}
$$
are linearly independent, 
the $B(n,r)$ are uniquely determined by the above equation. Hence
the Main Theorem is equivalent to:
$$
{{Z(n;2, \dots ,2,2+a; 0, \dots , 0)} \over
{Z(n;2, \dots ,2,2; 0, \dots , 0)}}=
{{q^{-a/2} } \over { {{3n-2} \choose {n-1}}}} 
\sum_{r=1}^{n} {{n+r-2} \choose {n-1}}{{2n-1-r} \choose {n-1}}
 [2+a]^{r-1}[1+a]^{n-r} \quad.
$$
 
By replacing $n$ by $n+1$, and changing the summation on $r$ to start
at $0$, we get that it suffices to prove:
$$
{{Z(n+1;2, \dots ,2,2+a; 0, \dots , 0)} \over
{Z(n+1;2, \dots ,2,2; 0, \dots , 0)}}=
{{q^{-a/2} } \over { {{3n+1} \choose {n}}}} 
\sum_{r=0}^{n} {{n+r} \choose {n}}{{2n-r} \choose {n}}
 [2+a]^{r}[1+a]^{n-r} \quad.
$$
 
Let $\tilde Z(n; x_1 , \dots , x_n; y_1 , \dots , y_n)$ denote
the right hand side of the Izergin-Korepin formula (Theorem 6 of [K]).
First replace $n$ by $n+1$.
Then, taking $x_i=2+i \eps$, for $i=1, \dots ,n$,
$x_{n+1}=2+a+(n+1)\eps$, and $y_j=-(j-1)\eps$, $j=1, \dots ,n+1$,
yields, after cancellation,
$$
{{\tilde Z(n+1\,;\,2+\eps, \dots ,2+n\eps\,,\,2+a+(n+1)\eps\,;\,
 0,-\eps, \dots , -n\eps)}
 \over
{\tilde Z(n+1\,;\,2+\eps, \dots ,2+n\eps\,,\,2+(n+1)\eps\,;\,
 0,-\eps, \dots , -n\eps)} }=
$$
$$
=
q^{-a/2} \prod_{j=0}^{n}
{{[2+a+(n+1+j)\eps] [1+a+(n+1+j)\eps] } 
\over
{[2+(n+1+j)\eps] [1+(n+1+j)\eps] } } \cdot
\prod_{j=0}^{n-1}
{{[ (n-j)\eps] } 
\over
{[ a+(n-j)\eps] }  } \cdot
{{det M_{n+1}(a)} \over {det M_{n+1}(0)}} \quad ,
$$
where $M_{n+1}(a)=(m_{i,j} \,,\, 0\leq i,j \leq n)$
is the $(n+1) \times (n+1)$ matrix, defined as follows:
$$
m_{i,j}= \cases{
1/\left ([2+(i+j+1)\eps][1+(i+j+1)\eps] \right), &  if $0 \leq i \leq n-1$,
                                               $0 \leq j \leq n$; \cr
1/\left([2+a+(n+j+1)\eps][1+a+(n+j+1)\eps]\right), &  if $i=n$,
                                               $0 \leq j \leq n$. \cr} 
$$
 
Taking the limit $\eps \rightarrow 0$, replacing $q^{\eps}$ by $s$,
$q^a$ by $X$, 
setting $w:=e^{\sqrt{-1} \pi/3}$,
evaluating the limit whenever possible, and cancelling
out whenever possible, reduces our task to proving the
following identity:
$$
\lim_{s \rightarrow 1} 
\left \{ {
 {(1-s)^n \det N_{n+1}(X)} 
     \over
 {\det N_{n+1}(1)}} \right\}
$$
$$
={
 {-(1-X)^n (\sqrt{-3})^{n+2} w^{-n}}
\over 
 {n!(1+X+X^2)^{n+1} {{3n+1} \choose {n}} }} \cdot
\sum_{r=0}^{n} w^{-r} {{n+r} \choose {n}}{{2n-r} \choose {n}}
(1+wX)^r(1-w^2 X)^{n-r} \quad .
\eqno(NotYetDone)
$$
 
Here the matrix $N_{n+1}(X)$ is  $M_{n+1}(a)$ divided by $3$, to wit:
$N_{n+1}(X)=(p_{i,j}, 0\leq i,j \leq n)$ is the 
$(n+1) \times (n+1)$ matrix, defined as follows. For the
first $n$ rows we have:
$$
p_{i,j}= {{1-s^{i+j+1}} \over {1-s^{3(i+j+1)}}}
\quad ,\quad (\, 0 \leq i \leq n-1 \,\,,\,\, 0 \leq j \leq n \,) \quad ,
$$
while for the last row we have:
$$
p_{n,j}={{1-X s^{n+j+1}} \over {1-X^3 s^{3(n+j+1)}}} \quad, \quad
(\,0 \leq j \leq n\,) \quad .
$$
We are left with the task of computing the determinant of $N_{n+1}(X)$, or
at least the limit on the left of $(NotYetDone)$.
\eject
{\bf A Short Course on Orthogonal Polynomials}
 
I will only cover what we need here. Of the many available
accounts, Ch. 2 of [Wilf] is especially recommended. For the
present purposes, section IV of [D] is most pertinent.
 
{\bf Theorem OP.}$-$ Let $T$ be any linear functional (`umbra') on
the set of polynomials, \break and let $c_i:=T(x^i)$ be its so-called
{\it moments}. Let 
$$
\Delta_n:=\,\det
\pmatrix{
c_0&c_1&\dots&\dots&c_n\cr
c_1&c_2&\dots&\dots&c_{n+1}\cr
\dots&\dots&\dots&\dots&\dots\cr
c_n&c_{n+1}&\dots&\dots&c_{2n}\cr
} \quad .
$$
If $\Delta_n \neq 0$, for $n \geq 0$, then there is a unique sequence
of {\it monic} polynomials $P_n(x)$,
where the degree of $P_n(x)$ is $n$, that are orthogonal with respect
to the functional $T$:
$$
T(P_n(x)P_m(x))= 0 \quad if \quad m \neq n \quad .
$$
Furthermore, these polynomials $P_n(x)$ are given `explicitly'
by:
$$
P_n(x)=
{{1} \over {\Delta_{n-1} }}
\det \pmatrix{
c_0&c_1&\dots&\dots&c_n\cr
c_1&c_2&\dots&\dots&c_{n+1}\cr
\dots&\dots&\dots&\dots&\dots\cr
c_{n-1}&c_{n}&\dots&\dots&c_{2n-1}\cr
1&x&x^2&\dots&x^{n}\cr
} \quad .
\eqno(General Formula)
$$
{\bf Proof:} \halmos.
 
{\bf Corollary 1.}$-$
$$
T(x^n P_n(x))={{\Delta_{n}} \over {\Delta_{n-1}} } \quad for \quad n \geq 1
 \quad .
$$
{\bf Proof:} \halmos.

{\bf Corollary 2.}$-$ If $S$ is another linear functional, 
$\,d_i:=S(x^i)\,$, and
$$
\Gamma_n:=
\det
\pmatrix{
c_0&c_1&\dots&\dots&c_n\cr
c_1&c_2&\dots&\dots&c_{n+1}\cr
\dots&\dots&\dots&\dots&\dots\cr
c_{n-1}&c_{n}&\dots&\dots&c_{2n-1}\cr
d_0&d_1&d_2&\dots&d_n\cr
} \quad .
$$
Then
$$
{{\Gamma_{n}} \over {\Delta_{n}} }={{S(P_n(x))} \over {T(x^n P_n(x))}}
 \quad for \quad n \geq 1
 \quad .
$$
{\bf Proof:} \halmos.
 
For a long time it was believed that theorem {\bf OP} was of only
theoretical interest, and that, given the moments, it was impractical
to actually find the polynomials $P_n(x)$, by evaluating the determinant.
This conventional wisdom was conveyed by Richard Askey, back in the
late seventies, to Jim Wilson, who was then studying under him.
Luckily, Wilson did not take this advice. Using
Theorem {\bf OP} lead him to beautiful
results[Wils], which later lead to the 
celebrated Askey-Wilson polynomials[AW]. Jim Wilson's
independence was later wholeheartedly endorsed by Dick Askey, who said:
{\it ``If an authority in the field tells you that
a certain approach is worth trying, listen to them. If they tell you
that a certain approach is {\bf not} worth trying, {\bf don't} listen
to them''}.
 
The `uselessness' of $(GeneralFormula)$
was still proclaimed a few years later, by
yet another authority, Jean Dieudonn\'e, who said ([D], p. 11):
{\it ``La formule g\'en\'erale [$(General Formula)$] donnant les $\dots$
sont impraticables pour le calcul explicite  $\dots$''}.
 
There is another way, by which $(General Formula)$, 
and its immediate corollaries {\bf 1} and {\bf 2} could be useful.
Suppose that we know, {\it by other means}, that a certain
set of explicitly given
monic orthogonal polynomials $Q_n(x)$ are orthogonal
with respect to the functional $T$, i.e. $T(Q_n Q_m)=0$ whenever $n \neq m$.
Then by uniqueness $Q_n=P_n$. If we are also able, using the
explicit expression for $Q_n(x)$, to find $T(x^n Q_n(x))$, then
Corollary {\bf 1} gives a way to {\it explicitly evaluate} the
Hankel determinant $\Delta_n$. If we are also able
to explicitly compute $S(Q_n(x))$, then we would be able
to evaluate the determinant $\Gamma_n$. This would be
our strategy in the evaluation of the determinants
on the left of $(NotYetDone)$, but we first need to digress again.
 
{\bf A Lean and Lively Course in q-Calculus}
 
Until further notice,
$$
(a)_n:=(1-a)(1-qa)(1-q^2 a) \dots (1-q^{n-1}a) \quad .
$$
 
If I had my way, I would ban $1-$Calculus from the Freshman
curriculum, and replace it by $q-$Calculus. Not only is it
more fun, it also describes nature more accurately. The
traditional calculus is based on the fictitious 
notion of the real line. It is now known that the universe
is quantized, and if you are at point $x$, then the
points that you can reach are in geometric progression
$q^i x$, in accordance with Hubble expansion. The true value
of $q$ is almost, but {\it not quite} $1$, and is a universal
constant, yet to be determined.
 
The $q$-derivative, $D_q$, is defined by
$$
D_q f(x) := {{f(x)-f(qx)} \over {(1-q)x}} \quad.
$$
The reader should verify that
$$
D_q x^a =  {{1-q^a} \over {1-q}} x^{a-1} \quad ,
$$
and the {\it product rule}:
$$
D_q [f(x) \cdot g(x)]= f(x) \cdot D_q g(x) + D_q f(x) \cdot g(qx) \quad .
\eqno(Product  Rule)
$$
 
The q-analog of integration, independently discovered by 
J. Thomae and the Rev. F.H. Jackson (see [AA], [GR]), is given by
$$
\int_{0}^{a} f(x) \,d_q x :=
a(1-q) \sum_{r=0}^{\infty} f(a q^r) q^r \quad ,
$$
and over a general interval:
$$
\int_{c}^{d} f(x) \,d_q x :=
\int_{0}^{d} f(x) \,d_q x -\int_{0}^{c} f(x) \,d_q x \quad .
$$
 
The reader is invited to use telescoping to prove the
{\bf Fundamental Theorem of q-Calculus}
$$
\int_{c}^{d} D_q F(x) \, d_q x \,= \,F(d)-F(c) \quad .
$$
 
Combining the Product Rule and the Fundamental Theorem, we
have:
 
{\bf q-Integration by Parts}: If $f(x)$ or $g(x)$ vanish at
the endpoints $c$ and $d$, then
$$
\int_{c}^{d} f(x)\cdot D_q g(x) \,d_q x = 
-\int_{c}^{d} D_q f(x) \cdot g(qx) \,d_q x \quad .
$$
 
{\bf Corollary:} If $g(q^i x)$ vanish at
the endpoints $c$ and $d$, for $i=0,1, \dots , n-1$, then
$$
\int_{c}^{d} f(x) \cdot D_q^n g(x) \,\,d_q x = 
(-1)^n \int_{c}^{d} D_q^n f(x) \cdot g(q^nx) \,\,d_q x \quad .
$$
 
Even those who still believe in $1$-Calculus can use
$q$-Calculus to advantage. All they have to do is let $q \rightarrow 1$
at the end.
 
The q-analog of $\int_{s}^{1} x^a dx = (1-s^{a+1})/(a+1)$ is
$$
{{1} \over {1-q}}
\int_{s}^{1} x^a d_q x = {{1-s^{a+1}} \over {1-q^{a+1}}} \quad .
\eqno(q-Moment)
$$
 
Now this looks familiar! Letting $q:=s^3$ in the definition of 
$N_{n+1}(X)$, given right after $(NotYetDone)$,
(this new `$q$' has nothing to do with the former Kuperberg $q$), we
see that the matrix $N_{n+1}(1)$ is the Hankel matrix of the
moments with respect to the functional
$$
T(f(x)):={{1} \over {1-q}} \int_{s}^{1} f(x) d_q x \quad.
$$
So all we need is to come up with orthogonal polynomials
with respect to the `q-Lebesgue-measure', over the
interval $[s,1]$.
 
{\bf q-Legendre Polynomials}
 
The ordinary Legendre polynomials, over an interval $(a,b)$ may
be defined in terms of the Rodrigues formula
$$
{{n!} \over {(2n)!}} D^n \{(x-a)^n (x-b)^n\} \quad .
$$
The orthogonality follows immediately by integration by parts.
This leads naturally to the $q$-analog,
$$
Q_n(x;a,b):=
{{(1-q)^n} \over {(q^{n+1})_n}}
 D_q^n \{(x-a)(x-qa) \dots (x-a q^{n-1})
\cdot (x-b)(x-qb) \dots (x-b q^{n-1})\} \quad .
$$
Using $q$-integration by parts repeatedly, it follows immediately that
the $Q_n(x;a,b)$ are orthogonal w.r.t. to $q-$integration over
$(a,b)$. The classical case $a=-1$, $b=1$ goes back to Markov.
Askey and Andrews[AA] generalized these to q-Jacobi polynomials,
and Askey and Wilson[AW] found the {\it ultimate} generalization.
While at present I don't see how to apply these more general
polynomials to combinatorial enumeration, I am sure that such
a use will be found in the future, and all enumerators are urged
to read [AA], [AW], and the modern classic [GR].
 
Going back to the determinant $N_{n+1}(X)$ of  $(NotYetDone)$, 
we also need to introduce the functional, defined on monomials by:
$$
S(x^j)={{1-Xs^{n+j+1}} \over {1-X^3 q^{n+j+1}}} \quad ,
$$
and extended linearly.
 
Let $X:=q^{\alpha/3}$. Then (recall that $s=q^{1/3}$):
$$
S(x^j)={{1-s^{\alpha+n+j+1}} \over {1- q^{n+j+\alpha+1}}}=
{{1} \over {1-q}} \int_{s}^{1} x^{\alpha+n} x^j \,d_q x \quad. 
$$
By linearity, for any polynomial $p(x)$:
$$
S(p(x))= {{1} \over {1-q}} \int_{s}^{1} x^{\alpha+n} p(x) \,d_q x \quad. 
$$
 
Using Corollary {\bf 2} of $(General Formula)$, we get
$$
{{ \det N_{n+1}(X)} \over  {\det N_{n+1}(1)}}=
{{\int_{s}^{1} x^{\alpha+n} P_n(x) \,d_q x } \over
{\int_{s}^{1} x^{n} P_n(x)\, d_q x } } \quad ,
\eqno(Almost Done)
$$
where $P_n(x)$ is now the $q-$Legendre polynomial over $[s,1]$, $Q_n(x;s,1)$
and $s=q^{1/3}$. In other words:
$$
P_n(x):=
{{(1-q)^n} \over { (q^{n+1})_n}} D_q^n \{(x-1)(x-q) \dots (x- q^{n-1}) \cdot
(x-s)(x-qs) \dots (x-s q^{n-1})\} \quad .
$$
\eject
{\bf DENOUEMENT}
 
It remains to compute the right side of $(Almost Done)$.
Let's first do the denominator.
 
{\bf Proposition Bottom.}$-$
$$
{{1} \over {1-q}}
\int_{s}^{1} x^n P_n(x) \,\,d_q x =
{{q^{n^2} (q)_n^2 (q^{-n}s)_{2n+1} } \over
{(q^{n+1})_n (q^{n+1})_{n+1} }} \quad.
$$
{\bf First Proof:} Use $q$-integration by parts,
$n$ times (i.e. use the above corollary). The resulting $q$-integral
is the famous $q$-Vandermonde-Chu sum,
that evaluates to the right side. See [GR], or use {\tt qEKHAD} accompanying
[PWZ]. \halmos
 
{\bf Remark:} Proposition Bottom, combined with Corollary 1 of
($General Formula$) gives an alternative evaluation of 
Kuperberg's determinant $N_n(1)$, needed in [K].
 
{\bf Second Proof:} Don't get off the shoulders of Greg Kuperberg yet.
Use his evaluation, and Corollary 1 of ($General Formula$).\halmos
 
{\bf Proposition Top.}$-$ Recalling that $X=q^{\alpha/3}$, we have
$$
{{1} \over {1-q}}
\int_{s}^{1} x^{n+\alpha} P_n(x) \,d_q x =
{{(-1)^n (qX^3)_n } \over {(q^{n+1})_n  }} \cdot
\sum_{k=0}^{\infty}
 q^k X^{3k} \prod_{r=0}^{n-1} (q^{n+k}-q^r) (q^{n+k}-q^{r+1/3})
$$
$$
-{{(-1)^n (qX^3)_n } \over {(q^{n+1})_n  }} \cdot \sum_{k=0}^{\infty} q^{k+1/3}
 X^{3k+1} \prod_{r=0}^{n-1} (q^{n+k+1/3}-q^r) (q^{n+k+1/3}-q^{r+1/3})
 \quad.
$$
{\bf Proof:} Let 
$$
F_n(x):=
{{(1-q)^n} \over { (q^{n+1})_n}} (x-1)(x-q) \dots (x- q^{n-1})
\cdot (x-s)(x-qs) \dots (x-s q^{n-1}) \quad ,
$$
so that $P_n(x)=D_q^n F_n(x)$. Since $F_n(q^ix)$ vanish
for $i=0, \dots ,n-1$ at both $x=1$ and $x=s$, we have by
$q$-integrating by parts $n$ times (the above corollary), that 
$$
\int_{s}^{1} x^{n+\alpha}\cdot P_n(x) \, d_q x =
\int_{s}^{1} x^{n+\alpha} \cdot D_q^n F_n(x) \, d_q x =
(-1)^n \int_{s}^{1} D_q^n \{x^{n+\alpha}\} \cdot F_n(q^n x) \,d_q x
$$
$$
={{(-1)^n(q^{\alpha+1})_n} \over {(1-q)^n}}
  \int_{s}^{1} x^{\alpha} F_n(q^n x) d_q x \quad .
$$
Now use the definition of $q$-integration over $[s,1]$ and
replace $q^{\alpha}$ by $X^3$, to complete the proof. \halmos
 
To compute the right side of $(Almost Done)$, we only need to
divide the expression given by Proposition Top by the
expression given by Proposition Bottom.
Doing this, multiplying by $(1-s)^n=(1-q^{1/3})^n$, and taking
the limit $q \rightarrow 1$, we get that the left side of
$(NotYetDone)$ is (Warning, now we are safely back in $1-$land,
so from now $(a)_n:=a(a+1) \dots (a+n-1)$, the ordinary
rising factorial):
$$
{{(-1)^n (1-X^3)^n (2n+1)!} \over {3^n n!^3 (-n+1/3)_{2n+1}}}
\cdot
\left( \sum_{k=0}^{\infty} (k+1)_n (k+2/3)_n X^{3k} -
\sum_{k=0}^{\infty} (k+1)_n (k+4/3)_n X^{3k+1} \right) \quad .
$$
 
After trivial cancellations, equation $(NotYetDone)$ boils down to
$$
({{1-X^3} \over {1-X}})^{2n+1}
{{ (-1)^n (3n+1)!} \over {3^{n+1} n!^3 (-n+1/3)_{2n+1}}}
\cdot
 \sum_{k=0}^{\infty} (k+1)_n (k+2/3)_n X^{3k} 
$$
$$
-({{1-X^3} \over {1-X}})^{2n+1}
{{ (-1)^n (3n+1)!} \over {3^{n+1} n!^3 (-n+1/3)_{2n+1}}}
\sum_{k=0}^{\infty} (k+1)_n (k+4/3)_n X^{3k+1} 
$$
$$
=
(\sqrt{-3})^n
\sum_{r=0}^{n} w^{-r-n} {{n+r} \choose {n}}{{2n-r} \choose {n}}
(1+wX)^r(1-w^2 X)^{n-r} \quad .
\eqno(Done)
$$
 
This was given to {\tt EKHAD}, the Maple package accompanying [PWZ].
{\tt EKHAD} found a certain linear homogeneous second order recurrence in $n$
that is satisfied by both sums on the left of $(Done)$ (and hence
by their difference), and also by the right side. It remains
to prove that both sides of $(Done)$ agree at $n=0,1$, 
which Maple did as well, even though it could be done by any human.
 
The input and output files are
obtainable by anonymous {\tt ftp} to 
{\tt ftp.math.temple.edu} (login as {\tt anonymous}), {\tt cd}-ing to directory
{\tt /pub/zeilberg/refined}, and {\tt get}ing 
files {\tt inDone} and {\tt outDone} respectively.
{\tt EKHAD} (and {\tt qEKHAD})
are available at directory {\tt /pub/zeilberg/programs}.
Alternatively, on the Web, go to my Home Page given at footnote $1$,
and click in the appropriate places. Of course, your computer should
be able to reproduce file {\tt outDone}. Once you have downloaded
{\tt EKHAD} and {\tt inDone} into a directory, type:
{\tt maple -q<inDone>outDone}. After $380$ seconds of CPU time,
{\tt outDone} would be ready.
 
{\bf References}
 
[AA] G.E. Andrews, and R. Askey, {\it Classical orthogonal polynomials},
in:``Polyn\^omes Orthogonaux et Applications'' (Proceedings, Bar-Le-Duc 1984),
edited by C. Brezinski et. al, Lecture Notes in Mathematics {\bf 1171}
36-62, Springer-Verlag, Berlin, 1985.
 
[AW] R. Askey and J. Wilson, {\it ``Some basic hypergeometric orthogonal
Polynomials that generalize Jacobi polynomials''}, Memoirs of the Amer.
Math. Soc. {\bf 319}, 1985.
 
[D] J. Dieudonn\'e, {\it Fractions continuees et polyn\^omes orthogonaux
dans l'oeuvre de E.N. Laguerre},
in:``Polyn\^omes Orthogonaux et Applications'' (Proceedings, Bar-Le-Duc 1984),
edited by C. Brezinski et. al, Lecture Notes in Mathematics {\bf 1171},
1-15, Springer-Verlag, Berlin, 1985.
 
[GR] G. Gasper and M. Rahman,{\it ``Basic Hypergeometric Series''},
Encyclopedia of Mathematics and its applications {\bf 35},
Cambridge University Press, Cambridge, England, 1990.
 
[KBI] V.E. Korepin, N.M. Bogoliubov and A.G. Izergin, 
{\it ``Quantum Inverse Scattering and Correlation Function''},
Cambridge University Press, Cambridge, England, 1993.
 
[K] Greg Kuperberg, {\it Another proof of the alternating sign matrix
conjecture}, Inter. Math. Res. Notes, to appear; {\tt AMS PPS \#199508-05-001};
Available from {\tt http://www.math.yale.edu/users/greg}.
 
[MRR1] W.H. Mills, D.P. Robbins, and H.Rumsey,
{\it Proof of the Macdonald conjecture} , Invent. Math. {\bf 66}(1982),
73-87.
 
[MRR2] W.H. Mills, D.P. Robbins, and H.Rumsey, {\it Alternating sign matrices
and descending plane partitions}, J. Combin. Theo. Ser. A, {\bf 34}(1983),
340-359.
 
[PWZ] M. Petkovsek, H. Wilf, and D. Zeilberger, {\it ``A=B''}, A.K. Peters,
Wellesley, 1996.
 
[S] R.P. Stanley {\it A baker's dozen of conjectures concerning
plane partitions} , in: COMBINATOIRE ENUMERATIVE, ed. by G. Labelle
and P. Leroux, Lecture Notes in Mathematics {\bf 1234}, Springer,
Berlin, 1986.
 
[Wilf] H. Wilf, {\it ``Mathematics for the Physical Sciences''},
Dover,  New York, 1978. Originally published by John Wiley, 1962.
 
[Wils] J. Wilson, {\it ``Hypergeometric series, recurrence relations,
and some new orthogonal functions''}, Ph.D. thesis, Univ. of Wisconsin,
Madison, 1978.
 
[Z] D. Zeilberger, {\it Proof of the alternating sign matrix conjecture},
Electronic J. of Combinatorics {\bf 3(2)}(1996)
[Foata Festschrift], R13. 
 
\bye